\documentclass[a4paper,12pt]{article}
\usepackage{amsmath}
\usepackage{amssymb}
\usepackage{amsfonts}
\begin{document}
\newtheorem{theoremenv}{Theorem}
\newtheorem{propositionenv}[theoremenv]{Proposition}
\newtheorem{lemmaenv}[theoremenv]{Lemma}
\newtheorem{corollaryenv}[theoremenv]{Corollary}
\newtheorem{claimenv}[theoremenv]{Claim}
\newtheorem{constructionenv}[theoremenv]{Construction}
\newtheorem{conjectureenv}[theoremenv]{Conjecture}
\newtheorem{questionenv}[theoremenv]{Question}
\newtheorem{definitionenv}[theoremenv]{Definition}
\newsavebox{\proofbox}
\savebox{\proofbox}{\begin{picture}(7,7)%
  \put(0,0){\framebox(7,7){}}\end{picture}}

\def\proof{\noindent\textbf{Proof. }}
\def\endproof{\hfill{\usebox{\proofbox}}}
\def\remark{\noindent\textbf{Remark. }}
\def\scs{\scriptstyle}
\def\boxeq{\tag*{\usebox{\proofbox}}}

\newcommand{\Hom}{\mbox{Hom}}
\newcommand{\Cov}{\mbox{Cov}}
\newcommand{\Var}{\mbox{Var}}
\newcommand{\Supp}{\mbox{Supp}}
\newcommand{\Diam}{\mbox{Diam}}
\newcommand{\Span}{\mbox{Span}}
\newcommand{\Prob}{\mbox{Prob}}
\newcommand{\RanProb}{\mathbb{P}_{\mbox{ran}}}
\newcommand{\SetProb}{\mathbb{P}_{\mbox{set}}}
\newcommand{\RegProb}{\mathbb{P}_{\mbox{reg}}}

    \begin{center}
    \Large{Counting sets with small sumset, and the clique number of random Cayley graphs}\\[15pt]
    \normalsize
     Ben Green\footnote{Supported by a grant from the Engineering and Physical Sciences Research Council of the UK and a Fellowship of Trinity College Cambridge.}\\[30pt]
    \end{center}
 \begin{abstract} \noindent Given a set $A \subseteq \mathbb{Z}/N\mathbb{Z}$ we may form a Cayley sum graph $G_A$ on vertex set $\mathbb{Z}/N\mathbb{Z}$ by joining $i$ to $j$ if and only if $i + j \in A$. We investigate the extent to which performing this construction with a random set $A$ simulates the generation of a random graph, proving that the clique number of $G_A$ is almost surely $O(\log N)$. This shows that Cayley sum graphs can furnish good examples of Ramsey graphs. To prove this result we must study the specific structure of set addition on $\mathbb{Z}/N\mathbb{Z}$. Indeed, we also show that the clique number of a random Cayley sum graph on $\Gamma = (\mathbb{Z}/2\mathbb{Z})^n$ is almost surely \textit{not} $O(\log |\Gamma|)$.
\end{abstract} 
\noindent\textbf{1. Introduction.} 
There are many papers in which a graph is constructed from a set in something like the following manner.
\begin{constructionenv}\label{cons1} Let $\Gamma$ be a finite abelian group with cardinality $N$ and let $A \subseteq \Gamma$. Define $G_A$ to be the graph on vertex set $\Gamma$ in which $ij \in E(G_A)$ precisely when $i + j \in A$.\end{constructionenv}
For any positive integer $m$ write $\mathbb{Z}_m = \mathbb{Z}/m\mathbb{Z}$. In this paper we will be interested in the two particular groups $\Gamma = \mathbb{Z}_N$ and $\Gamma = \mathbb{Z}_2^n$.\\[11pt]
\noindent\textbf{Example.} Suppose that $N$ is prime and let $A$ be the set of quadratic residues modulo $N$. The graph $G_A$ is then called the \textit{Paley sum graph} of order $N$, and is denoted by $P_N$. It is generally thought that the Paley graphs behave much like random graphs from the binomial model $G(N,1/2)$, though in many respects this statement is far from being proven. We remind the reader that in this model all (labelled) graphs on $N$ vertices occur with equal probability. Equivalently, a graph from $G(N,1/2)$ may be generated by picking each possible edge $ij$ independently at random with probability $1/2$. See \cite{Boll} for the basic definitions of random graph theory (and much more).\\[11pt]
In the light of the above example it seems natural to address the following question. 
\begin{questionenv}\label{mainq} To what extent does choosing a random subset $A \subseteq \Gamma$ and forming $G_A$ simulate the selection of a random graph on $N$ vertices? 
\end{questionenv} We call a graph generated using this model a \textit{random Cayley sum graph}.
The reader may have already noticed that there is one sense in which a random Cayley sum graph is definitely \textit{not} like a random graph -- it is always regular. However in a recent paper of Krivelevich, Sudakov, Vu and Wormald \cite{KSVW} a good quantity of evidence was collected to support the claim that a random regular graph of degree $N/2$ (say) behaves in many ways almost exactly like a random graph from the binomial model $G(N,1/2)$. \\[11pt]
In this paper we will investigate just one aspect of Question \ref{mainq}, the clique number of a random Cayley sum graph. There are clearly other ways in which we might compare a random Cayley sum graph to a random graph. An investigation of the second largest eigenvalue, for example, leads naturally to some classical questions on random trigonometric polynomials. We discuss this matter in a separate paper \cite{Green}.\\[11pt]
This is really a paper in combinatorial number theory. To see why, let $A \subseteq \Gamma$ and consider what it means for $X \subseteq \Gamma$ to span a clique in $G_A$. We must have $i + j \in A$ whenever $i,j$ are distinct elements of $X$. Now the set $\{i + j : i,j \in X, i \neq j\}$ is called the restricted sumset of $X$ with itself and is normally written $X \hat{+} X$. Thus we see that $X$ spans a clique in $G_A$ if and only if $X \hat{+} X \subseteq A$.\\[11pt]
 Now if $X$ is highly structured, for example if it is an arithmetic progression, then $|X \hat{+} X|$ is rather small compared to $|X|$ and the probability that $G_A$ contains the clique spanned by $X$ will be quite large. This does not mean, however, that the clique number of a random Cayley graph is very large since there are rather few arithmetic progressions. \\[11pt]
Let $k,m \in \mathbb{N}$, and make the following definition.
\begin{definitionenv}
$S_k^m(\Gamma)$ is the collection of sets $X \subseteq \Gamma$ with $|X| = k$ and $|X \hat{+} X| = m$. 
\end{definitionenv}
Observe that $S_k^m(\Gamma)$ is empty if $m < k-1$. Now one can easily write down a formula for $E(\Gamma,k)$, the expected number of $k$-cliques in a random Cayley graph on $\Gamma$. Indeed one has
\begin{eqnarray*} E(\Gamma;k) & = & \sum_{X \subseteq \Gamma : |X| = k} \mathbb{P}(\mbox{$X$ is a clique in $G_A$}) \\ & = & \sum_{m \geq k-1}\sum_{X \in S_k^m(\Gamma)} \mathbb{P}(X \hat{+} X \subseteq A) \\ & = & \sum_{m \geq k-1} \sum_{X \in S_k^m(\Gamma)} 2^{-m} \\ & = & \sum_{m \geq k-1} |S_k^m(\Gamma)|2^{-m}.\end{eqnarray*}
Now suppose we could show, for some function $F = F_{\Gamma} : \mathbb{N} \rightarrow \mathbb{N}$, that $E(\Gamma,F(N)) = o(1)$. Writing $\omega$ for the clique number and $\Omega_k$ for the number of cliques of size $k$ we would then have
\[ \mathbb{P}(\omega \geq F(N)) \; = \; \mathbb{P}(\Omega_{F(N)} \geq 1) \; \leq \; E(\Gamma,F(N)) \; = \; o(1)\] by Markov's inequality. That is, the clique number of a random Cayley graph on $\Gamma$ is almost surely at most $F(N)$. As a result of these observations our attention is drawn to the problem of finding effective upper bounds for $|S_k^m(\Gamma)|$.\\[11pt]
\noindent\textbf{Further notation and background.} Let $|\Gamma| = N$ and consider the sample space $\Omega_N$ consisting of all graphs on $N$ vertices. If $E = E_N$ is an event in this sample space then we write 
\begin{itemize}
\item $\RanProb^{(N)}(E)$ for the probability that a random graph from $G(N,1/2)$ satisfies $E$;
\item $\RegProb^{(N)}(E)$ for the probability that a graph selected from $\mathcal{G}_{N,\lfloor N/2\rfloor}$ satisfies $E$, where $\mathcal{G}_{N,\lfloor N/2\rfloor}$ puts the uniform distribution on all the $\lfloor N/2 \rfloor$-regular graphs on $N$ labelled vertices;
\item  $\mathbb{P}_{\Gamma}(E)$ for the probability that if $A \subseteq \Gamma$ is chosen randomly,  then the Cayley sum graph $G_A$ satisfies $E$.
\end{itemize}
The following well-known result concerning the clique number of a random graph (see, for example, \cite{Boll}) sets our later results in context. In this paper, all logs are to the base 2 unless otherwise specified.
\begin{theoremenv} For any fixed $\epsilon > 0$ we have
\[\lim_{N \rightarrow \infty}\mathbb{P}^{(N)}_{\mbox{\emph{ran}}}\left((2 - \epsilon)\log N \; \leq \; \omega \; \leq \; (2 + \epsilon)\log N\right) \; = \; 1.\]
\end{theoremenv}
In \cite{KSVW} the same result is proved for random regular graphs:
\begin{theoremenv}
For any fixed $\epsilon > 0$ we have
\[\lim_{N \rightarrow \infty}\mathbb{P}^{(N)}_{\mbox{\emph{reg}}}\left((2 - \epsilon)\log N \; \leq \; \omega \; \leq \; (2 + \epsilon)\log N\right) \; = \; 1.\]
\end{theoremenv}
Questions equivalent to estimating the clique number of a random Cayley sum graph appear in various places in the literature\footnote{The author would like to thank Benny Sudakov for drawing his attention to the two papers \cite{AAAS,AO}.}, such as \cite{AAAS,Alon,AO,Choi,Ruzsa4}. To the best of my knowledge the following result, which is essentially given in \cite{AAAS}, is the best known. 
\begin{theoremenv}\label{thm4} 
There is an absolute constant $C$ such that \[ \lim_{N \rightarrow \infty}\mathbb{P}_{\mathbb{Z}_N}\left(\omega \; \leq \; C(\log N)^2\right) \; = \; 1. \]
\end{theoremenv}
A similar result is mentioned in \cite{AO}, where it is remarked that it would be interesting to try and improve on it, and also, without proof, in \cite{Ruzsa4}.\\[11pt] 
Let us conclude this section by making a remark concerning our terminology. We use the term \textit{Cayley sum graph} to make it clear that the objects we are considering are very slightly different from what are normally called Cayley graphs. In an abelian group $\Gamma$ a (true) Cayley graph is formed as follows. Take a symmetric set $A$, that is to say a set in which $a \in A$ implies that $-a \in A$. Form a graph on vertex set $\Gamma$ by joining $i$ to $j$ if and only if $i - j \in A$. The group $\Gamma$ acts on the resulting graph in an obvious way, and for this reason true Cayley graphs are perhaps more natural than Cayley sum graphs. From the rather number-theoretical point of view of this paper, however, the treatment of Cayley sum graphs seems to be a little tidier. The methods of this paper apply to true Cayley graphs with minimal modifications. When $\Gamma = \mathbb{Z}_2^n$, of course, the two types of graph actually coincide.\\[11pt] 
 \noindent\textbf{2. Statement of results.}
The main result of this paper is the following improvement of Theorem \ref{thm4}.
\begin{theoremenv}\label{mainthm1}
We have 
\[ \lim_{N \rightarrow \infty}\mathbb{P}_{\mathbb{Z}_N}\left(\omega \; \leq \; 160\log N\right) = 1.\]
\end{theoremenv}
An immediate consequence of Theorem \ref{mainthm1} is the following result in Ramsey theory, for which we do not know an independent proof:
\begin{theoremenv}
For all sufficiently large integers $N$ there exists a set $A \subseteq \mathbb{Z}_N$ for which the Cayley sum graph $G_A$ has no cliques or independent sets of size $160\log N$.
\end{theoremenv}
The behaviour of random Cayley sum graphs on $\mathbb{Z}_2^n$ is different, reflecting the fact that this group has much richer subgroup structure. We will prove the following.
\begin{theoremenv}\label{mainthm2}
There are absolute constants $C_1,C_2 > 0$ such that
\[ \lim_{N \rightarrow \infty}\mathbb{P}_{\mathbb{Z}_2^n}\left(C_1\log N \log\log N \; \leq \; \omega \; \leq \; C_2\log N \log\log N\right) \; = \; 1,\] where $N = 2^n = |\mathbb{Z}_2^n|$.
\end{theoremenv}
We make no effort to optimise the values of $C_1$ and $C_2$ in this result.\\[11pt]
\noindent\textbf{3. Freiman isomorphism classes.} A vital notion in this paper is that of a \textit{Freiman isomorphism}. Let $s$ be a positive integer, let $A$ and $B$ be subsets of (possibly different) abelian groups and let $\phi : A \rightarrow B$ be a map. Then we say that $\phi$ is a Freiman $s$-homomorphism if whenever $a_1,\dots,a_s,a_1',\dots,a_s' \in A$ satisfy
\begin{equation}\label{eq271} a_1 + a_2 + \dots + a_s \; = \; a_1' + a'_2 + \dots + a_s'\end{equation}
we have
\begin{equation}\label{eq272} \phi(a_1) + \phi(a_2) + \dots + \phi(a_s) \; = \; \phi(a_1') + \phi(a_2') + \dots + \phi(a_s').\end{equation}
If $\phi$ has an inverse which is also a $s$-homomorphism then we say that it is a Freiman $s$-isomorphism, and write $A \cong_s B$. This is the case precisely if \eqref{eq271} holds for any $a_i,a'_i$ satisfying \eqref{eq272}.\\[11pt]
This is a standard and very convenient terminology in additive number theory. One must, however, be slightly careful: we invite the reader to consider the obvious map $\phi : \{0,1\} \rightarrow \mathbb{Z}_2$, which is a one to one Freiman 2-homomorphism but not a Freiman isomorphism.\\[11pt]
In this paper we will mostly, but not exclusively, be concerned with Freiman $2$-homomorphisms. From now on we will refer to these simply as homomorphisms.\\[11pt]
In order to discuss the theory of Freiman homomorphisms as cleanly as possible we introduce some further definitions. These are very much in the spirit of \cite{LK}, and indeed Lemma \ref{lem12} below is a crude (though slightly more general) version of one of the results of that paper.\\[11pt]
Let $F$ be a field, let $k$ be a positive integer and let $W$ be a vector space over $F$. Let $V_k$ be a $k$-dimensional vector space over $F$ with basis $\{e_1,\dots,e_k\}$. For any sequence $\sigma = (a_1,\dots,a_k)$ of $k$ elements of $W$ we can define a linear map $\phi_{\sigma} : V_k \rightarrow W$ by
\[ \phi_{\sigma}\left(\sum \lambda_i e_i\right) \; = \; \sum \lambda_i a_i.\]
If $v \in \ker \phi_{\sigma}$ then we say that $v$ is \textit{satisfied by $\sigma$}. Let $s \geq 2$ be an integer. Call any vector of the form
\[ e_{i_1} + \dots + e_{i_s} - e_{j_1} - \dots - e_{j_s},\] where the $i$s and $j$s need not be distinct, an $s$-\textit{relation}.  Write $\mathcal{R}_s$ for the set of all $s$-relations in $V_k$. If $S \subseteq V_k$ write $\Span(S)$ for the $F$-linear span of the elements of $S$.\\[11pt]
Now for any subset $S$ and any subspace $U$ of a vector space one has $S \cap \Span(S \cap U) = S \cap U$. Thus in particular
\begin{equation}\label{eq273} \mathcal{R}_s \cap \Span\left(\mathcal{R}_s \cap \ker \phi_{\sigma}\right) \; = \; \mathcal{R}_s \cap \ker \phi_{\sigma}.\end{equation} That is to say, the set $\mathcal{R}_s \cap \ker \phi_{\sigma}$ of $s$-relations satisfied by $\sigma$ can be completely recovered from the subspace spanned by it.
\begin{lemmaenv}\label{lem2}
Suppose that $\sigma = (a_i)_{i=1}^k$ and $\sigma' = (a_i')_{i=1}^k$ are two sequences of distinct elements of $W$, and suppose that $\mbox{\emph{Span}}\left(\mathcal{R}_s \cap \ker \phi_{\sigma}\right) = \mbox{\emph{Span}}\left(\mathcal{R}_s \cap \ker \phi_{\sigma'}\right)$. Then the sets $A = \{a_1,\dots,a_k\}$ and $A' = \{a_1,\dots,a_k\}$ are $s$-isomorphic.
\end{lemmaenv}
\proof By \eqref{eq273} we have $\mathcal{R}_s \cap \ker \phi_{\sigma} = \mathcal{R}_s \cap \ker \phi_{\sigma'}$. That is, the set of $s$-relations satisfied by $\sigma$ and $\sigma'$ are the same. An $s$-isomorphism between $A$ and $A'$ is therefore given by the map $a_i \mapsto a_i'$.\endproof
\begin{lemmaenv}\label{lem12}
The number of $s$-isomorphism classes of subsets of $W$ of size $k$ is at most $k^{2sk}$.
\end{lemmaenv}
\proof Given any set $A = \{a_1,\dots,a_k\} \subseteq W$ we may form a sequence $\sigma = (a_1,\dots,a_k)$. By Lemma \ref{lem2} the number of $s$-isomorphism classes of sets $A$ is at most the number of subspaces of $V_k$ spanned by $s$-relations. There are at most $k^{2s}$ of these $s$-relations, and any such subspace is spanned by some subset of at most $k$ of them. The required bound follows.\endproof\\[11pt]
Now if $W,W'$ are two vector spaces over $F$ and if $A \subseteq W$ then the space $\Hom_s(A,W')$ of Freiman $s$-homomorphisms $\phi : A \rightarrow W'$ is a vector space over $F$. We call $r_F(A) = \dim \Hom_2(A,F) - 1$ the \textit{Freiman dimension} of $A$. It is easy to see that $\dim \Hom_2(A,W') = (r_F(A) + 1)^{\dim W'}$. In a short while we will prove the main result of this section, Lemma \ref{lem4}. Let us first prepare the ground with the following lemma from linear algebra.
\begin{lemmaenv}\label{lem4pre} Let $F$ be a field, let $A$ be a finite non-empty set, and let $\Phi$ be a subspace of $\Psi$, the vector space of all maps from $A$ to $F$. Write $d = \dim \Phi$. Then there exists a subset $\{a_1,\dots,a_d\} \subseteq A$ with the following properties.\\
\emph{(i)} For any $f_1,\dots,f_d \in F$ there is a unique $\phi \in \Phi$ with $\phi(a_i) = f_i$ \emph{(}$i = 1,\dots,d$\emph{)}.\\
\emph{(ii)} If $a \in A$ then there are $\lambda_1,\dots,\lambda_d \in F$ \emph{(}depending only on $a$\emph{)} such that $\phi(a) = \lambda_1\phi(a_1) + \dots + \lambda_d\phi(a_d)$ for any $\phi \in \Phi$. If the constant map, which sends everything in $A$ to $1_F$, lies in $\Phi$ then we may assume that $\lambda_1 + \dots + \lambda_d = 1$.
\end{lemmaenv}
\proof Let $A = \{a_1,\dots,a_k\}$. Choose as a basis for $\Psi$ the set $\{\alpha_1,\dots,\alpha_k\}$ of maps defined by the conditions $\alpha_i(a_j) = \delta_{ij}$. Now we may find some $k - d$ of the $\alpha_i$ which, when taken together with $\Phi$, span the whole of $\Psi$. Relabeling if necessary we may assume that these elements are $\alpha_{d+1},\dots,\alpha_k$, so that
\begin{equation}\label{eq771} \Psi \; = \; \Phi \oplus \Span(\alpha_{d+1},\dots,\alpha_k).\end{equation}
It follows that for any $f_1,\dots,f_d \in F$ we may find $\phi \in \Phi$ and $\mu_{d+1},\dots,\mu_k \in F$ so that 
\[ f_1\alpha_1 + \dots + f_d\alpha_d \; = \; \phi + \mu_{d+1}\alpha_{d+1} + \dots + \mu_k\alpha_k,\] implying that $\phi(a_i) = f_i$ ($i = 1,\dots,d$). To see that $\phi$ is unique, suppose that $\phi'(a_i) = \phi''(a_i)$ for some $\phi',\phi'' \in \Phi$ and all $i = 1,\dots,d$. Then $\phi' - \phi'' \in \Span(\alpha_{d+1},\dots,\alpha_k)$ and thus $\phi' - \phi'' = 0$ since the sum in \eqref{eq771} is direct.

As for assertion (ii), observe that any $x \in A$ induces a linear functional $t_x : \Phi \rightarrow F$ via the rule $t_x(\phi) = \phi(x)$. For a fixed $a \in A$ the $d+1$ linear functionals $t_a,t_{a_1},\dots,t_{a_d}$ are linearly dependent and therefore we can find $\mu,\mu_1,\dots,\mu_d \in F$, not all zero, so that \begin{equation}\label{eq980} \mu \phi(a) + \mu_1\phi(a_1) + \dots + \mu_d\phi(a_d) \; = \; 0\end{equation} for any $\phi \in \Phi$. If $\mu$ were equal to zero then we could easily contradict (i) by setting $f_i = 1$ and $f_j = 0$ for $j \neq i$. Therefore this is not the case and we may divide \eqref{eq980} through by $\mu$ to get (ii) (with $\lambda_i = -\mu^{-1}\mu_i$ for $i = 1,\dots,d$). If the constant map lies in $\Phi$ then \eqref{eq980} must hold for it, which implies that $\lambda_1 + \dots + \lambda_d = 1$.
\endproof
\begin{lemmaenv}\label{lem4} Let $W$ and $W'$ be finite-dimensional vector spaces over a field $F$ and let $A \subseteq W$ be a finite non-empty set of Freiman dimension $r = r_F(A)$. Then there is a subset $\{a_1,\dots,a_{r+1}\} \subseteq A$ with the following properties.\\
\emph{(i)} For any $w'_1,\dots,w'_{r+1} \in W'$ there is a unique $\phi \in \mbox{\emph{Hom}}_2(A,W')$ with $\phi(a_i) = w'_i$ \emph{(}$i = 1,\dots,r+1$\emph{)}. Moreover, if $w'_1,\dots,w_{r+1}'$ are linearly independent then $\phi$ is an isomorphism.\\
\emph{(ii)} If $a \in A$ then there are $\lambda_1,\dots,\lambda_{r+1} \in F$ \emph{(}depending only on $a$\emph{)} with $\lambda_1 + \dots + \lambda_{r+1} = 1_F$ and such that $\phi(a) = \lambda_1\phi(a_1) + \dots + \lambda_{r+1}\phi(a_{r+1})$ for any $\phi \in \mbox{\emph{Hom}}_2(A,W')$.\end{lemmaenv}
\proof Part (ii) and the first part of (i) are swift consequences of Lemma \ref{lem4pre}, (i). Although that lemma deals only with homomorphisms from $A$ to $F$, the required extension follows easily using the isomorphism  $\Hom_2(A,F^n) \cong (\Hom_2(A,F))^n$.

To prove the second assertion of (i), suppose that suppose that $\phi \in \Hom_2(A,W')$ is not a Freiman isomorphism. Then there are elements $a,a',a'',a''' \in A$ such that $\phi(a) + \phi(a') = \phi(a'') + \phi(a''')$ but $a + a' \neq a'' + a'''$. By part (ii) there are $u_i,u'_i,u''_i,u'''_i \in F$ such that $\psi(a) = \sum_{i = 1}^{r+1}u_i\psi(a_i)$ for all $\psi \in \Hom_2(A,V)$, for any vector space $V$ over $F$ (and similarly for $a',a'',a'''$). One has
\begin{equation}\label{eq377} \psi(a) + \psi(a') - \psi(a'') - \psi(a''') \; = \; \sum_{i = 1}^{r+1} (u_i + u'_i - u''_i - u'''_i)\psi(a_i)\end{equation} for all $\psi$.
Taking $V = W$ and $\psi = \mbox{id}$ we see that not all of the $u_i + u'_i - u''_i - u'''_i$ are zero. Therefore, setting $V = W'$ and $\psi = \phi$, one sees that
 $\phi(a_1),\dots,\phi(a_{r+1})$ must be linearly dependent. 
\endproof
\begin{corollaryenv}\label{cor1} Let $A \subseteq W$. Then $r = r_{F}(A)$ is the largest integer $d$ such that $A$ is isomorphic to a subset of $F^{d}$ which is not contained in an affine subspace. In particular if $A$ is a set of integers then $r_{\mathbb{Q}}(A) = d$ really is equal to dimension in the original sense of Freiman \emph{\cite{Fr}}.
\end{corollaryenv}
\proof Apply Lemma \ref{lem4} (ii) with $W' = W$ and $\phi = \mbox{id}$. It follows that there are $\{a_1,a_2,\dots,a_{r+1}\}$ such that $A$ is a subset of $\{\lambda_1a_1 + \dots + \lambda_{r+1}a_{r+1} : \lambda_1 + \dots + \lambda_{r+1} = 1_F\}$, an affine space over $F$ of dimension $r$. As $r_F(A)$ is clearly an isomorphism invariant, the same must be true of any set isomorphic to $A$ and it therefore follows that $d \leq r$. To see that $r \leq d$, pick linearly independent vectors $w'_1,\dots,w'_{r+1}$ in some space $W'$, and observe that $A$ is isomorphic to some set containing all of the $w'_i$ by Lemma \ref{lem4} (i).\endproof\\[11pt] 
\noindent\textbf{4. Isomorphism classes with small doubling.} As in the last section, let $W$ be a vector space over a field $F$.
The main result of this section is a bound for the number of mutually non-isomorphic sets $A \subseteq W$ of cardinality $k$ with small (restricted) doubling property. The cardinality of $A \hat{+} A$ is an isomorphism invariant and we say that an isomorphism class of sets has small doubling property if $m = |A \hat{+} A|$ is relatively small in comparison to $|A|$ for each set in that class. \\[11pt]
If $B \subseteq W$ and if $t$ is a non-negative integer than a $t$-\textit{extension} of $B$ is simply any set $B \cup E$ where $E = \{x_1,\dots,x_t\}$ and the $x_i$ are distinct elements of $W \setminus B$. 
\begin{propositionenv}\label{prop2}
Suppose that $A \subseteq W$ has sufficiently large cardinality $k = |A|$ and let $m = |A \hat{+} A|$. Then there exist $a^{\ast} \in A$ and $A_0,A_1 \subseteq A$ such that $|A_0| \leq 4k^{-1/15}m$, $|A\setminus A_1| \leq 4k^{-4/5}m$ and $a^{\ast} + A_1 \subseteq A_0 \hat{+} A_0$.
\end{propositionenv}
\proof Let $Q = \lfloor k^{1/5}\rfloor$ and $q = k^{-1/15}$. Draw a graph on vertex set $A$ by joining $x$ to $y$ if $x + y$ is $Q$-\textit{unpopular}, that is if the number of pairs $(w,z) \in A^2$, $w \neq z$, with $w + z = x + y$ is less than $Q$. The total number of edges in the graph is no more than $Qm$, and so some vertex $a^{\ast}$ has degree at most $2Qm/k$. This means that with at most $2Qm/k$ exceptions the sums $a^{\ast} + a$, $a \in A$ are $Q$-popular. Pick a set $X \subseteq A$ by choosing each $a \in A$ to lie in $X$ independently and at random with probability $q$. Let $Z = \{a \in A | a^{\ast} + a \notin X \hat{+} X\}$. If $s = a^{\ast} + a$ is $Q$-popular then there are at least $Q/3$ disjoint pairs $(b,c) \in A^2$ with $b + c = s$ and $b \neq c$. For each the probability that both $b$ and $c$ lie in $X$ is $q^2$, and so 
\[
\mathbb{P}(a \in Z) \; \leq \; (1 - q^2)^{Q/3} \; \leq \; e^{-q^2Q/3}.\]
It follows that 
\[ \mathbb{E}|Z| \; \leq \; \frac{2Qm}{k} + e^{-q^2Q/3}k \; \leq \; 3k^{-4/5}m.\]
Clearly we have $\mathbb{E}|X| = qk = k^{14/15}$ and so 
\begin{equation}\label{eq78} \mathbb{E}\left(|X| + k^{11/15}|Z|\right) \; \leq \; 4k^{-1/15}m.\end{equation}
Pick a specific $X \subseteq A$ for which $|X| + k^{11/15}|Z| \leq  4k^{-1/15}m$, and set $A_0 = X$, $A_1 = A \setminus Z$. Then $|A_0| \leq 4k^{-1/15}m$, $|A \setminus A_1| \leq 4k^{-4/5}m$ and, by construction, $a^{\ast} + A_1 \subseteq A_0 \hat{+} A_0$.\endproof
\begin{lemmaenv}\label{extlemma} Fix a non-negative integer $t$ and a set $B \subseteq W$ of cardinality $|B| = l$. Then the number of mutually non-isomorphic sets $A$ of cardinality $|A| = l + t$, such that there exists a subset $A_0 \subseteq A$ satisfying $A_0 \cong_3 B$, is at most $l^{3t^4}$.
\end{lemmaenv}
\proof Write $B = \{b_1,\dots,b_{l}\}$. Given $A_0 \subseteq A$ with the required properties write $A_0 = \{a_1,\dots,a_l\}$ and $A = A_0 \cup \{a_{l+1},\dots,a_{l+t}\}$ in a manner that respects the 3-isomorphism between $A_0$ and $B$; that is, for $i_1,i_2,i_3,i_4,i_5,i_6 \in \{1,\dots,l\}$
\[ b_{i_1} + b_{i_2} + b_{i_3} - b_{i_4} - b_{i_5} - b_{i_6} \; = \; 0\] if and only if
\[ a_{i_1} + a_{i_2} + a_{i_3} - a_{i_4} - a_{i_5} - a_{i_6} \; = \; 0.\]
The isomorphism class of $A$ is determined by the 2-relations in $V_{l+t}$ satisfied by $A$. With each possible relation we can associate its projection $w \in V_t$ onto the last $t$ coordinates of $V_{l+t}$. There are at most $t^4$ non-zero projections $w$; fix one of them, say $w = (1,0,\dots,0)$, and consider relations corresponding to this specific projection. If $A$ satisfies such a relation, say 
\begin{equation}\label{eq365}
a_{i_1} + a_{l+1} \; = \; a_{i_2} + a_{i_3}
\end{equation}
with $i_1,i_2,i_3 \in \{1,\dots,l\}$, then all other such relations are uniquely determined: they correspond to the equalities $a_{i_2} + a_{i_3} - a_{i_1} = a_{j_2} + a_{j_3} - a_{j_1}$ with $j_1,j_2,j_3 \in \{1,\dots,l\}$, hence to the equalities $b_{i_2} + b_{i_3} - b_{i_1} = b_{j_2} + b_{j_3} - b_{j_1}$ with $j_1,j_2,j_3 \in \{1,\dots,l\}$. Since there are at most $l^3$ choices for the relation \eqref{eq365}, there are at most $l^3$ choices for the whole system of relations corresponding to $w$ which $A$ satisfies. Moreover, it is not difficult to see that this conclusion is actually true for any possible $w \neq 0$. The result follows on observing that to determine the system of relations in $V_{l+t}$, satisfied by $A$, we must for each non-zero $w \in V_t$ determine the system of relations corresponding to $w$.\endproof\\[11pt]
The next lemma is an easy exercise.
\begin{lemmaenv}\label{lem36}
Suppose that $X \cong_6 X'$. Then $X \hat{+} X \cong_3 X'\hat{+}X'$ and any $B \subseteq X \hat{+} X$ is $3$-isomorphic to a subset of $X' \hat{+} X'$.
\end{lemmaenv}
\begin{propositionenv}\label{prop18a} Let $k$ be a sufficiently large positive integer and let $m \leq k^{31/30}$. Then the number of isomorphism classes of $A \subseteq F^k$ with $|A| = k$ and $|A \hat{+} A| \leq m$ is at most $(em/k)^k\exp\left(k^{31/32}\right)$.
\end{propositionenv}
\proof By Proposition \ref{prop2}, to any $A \subseteq F^k$ with $|A| = k$ and $|A \hat{+} A| \leq m$ there correspond $a^{\ast} \in A$ and $A_0,A_1 \subseteq A$ such that $|A_0| \leq 4k^{-1/15}m$, $|A_1| \geq |A| - 4k^{-4/5}m$ and $a^{\ast} + A_1 \subseteq A_0 \hat{+} A_0$. By Lemma \ref{lem12}, the number of possible 6-isomorphism classes to which $A_0$ can belong is $M < \exp\left(k^{30/31}\right)$; let $X_1,\dots,X_M$ be representatives of these classes. By Lemma \ref{lem36}, the set $a^{\ast} + A_1$ is 3-isomorphic to a non-empty subset of some $X_i\hat{+} X_i$ of cardinality at most $k$; consequently, the number of 3-isomorphism classes to which $A_1$ can belong is at most $Mk(em/k)^k$. By Lemma \ref{extlemma}, each of these 3-isomorphism classes gives rise, by $t$-extensions for any fixed $t \leq 4k^{-4/5}m$, to at most $\exp(k^{29/30})$ isomorphism classes. Altogether, we get at most
\[ \exp(k^{30/31}) \cdot k(em/k)^k \cdot 4k^{-4/5}m \cdot \exp(k^{29/30}) \; < \; (em/k)^k \exp(k^{31/32})\] possible isomorphism classes for $A$.\endproof\\[11pt]
\noindent\textbf{5. Avoiding very small doubling.} The bound of Proposition \ref{prop18a} becomes rather ineffective when $m \leq 7k$, say. This section contains a trick which allows us to avoid ever considering sets with extremely small doubling (see, however, the section entitled Concluding Remarks). The result we will actually use, Corollary \ref{trickcor}, is a consequence of the following more general proposition.
\begin{propositionenv}\label{trickprop} Let $k$ be a sufficiently large positive integer and let $A$ be a subset of an abelian group with $|A| = k$. Let $\epsilon > 0$ be sufficiently small. Then there is a set $B \subseteq A$ with 
\[ |B| \; \leq \; \frac{3\log(1/\epsilon)}{\epsilon}\sqrt{k}\] and $|B \hat{+} B| \geq (1 - \epsilon)k$.
\end{propositionenv}
\proof Let $p = 3\epsilon^{-1}\sqrt{\log(6/\epsilon)/2k}$. Choose a set $B \subseteq A$ by picking each element of $k$ independently and at random with probability $p$. By some rather standard tail estimates (see, for example, \cite[Theorem A.1.4]{AlonSpencer}) we have
\begin{equation}\label{eq717}\mathbb{P} \left(\left| |B| - pk \right| \; \geq \; \epsilon p k/3\right) \; < \; 2e^{-2\epsilon^2p^2k/9} \; = \; \epsilon/3.\end{equation}
Now for any $x \in B \hat{+} B$ write $s(x)$ for the number of ordered pairs $(b_1,b_2) \in B^2$ with $b_1 \neq b_2$ and $b_1 + b_2 = x$. The sum $M(B) = \sum_x s(x)^2$ equals the number of quadruples $(b_1,b_2,b_3,b_4) \in B^4$ with $b_1 + b_2 = b_3 + b_4$ and $b_1 \neq b_2$, $b_3 \neq b_4$. We call these summing quadruples.

Let us divide the summing quadruples $(a_1,a_2,a_3,a_4) \in A^4$ into two classes: those in which $a_1,a_2,a_3,a_4$ are all distinct and those in which either $a_1 = a_3, a_2 = a_4$ or $a_1 = a_4, a_2 = a_3$. Call this latter type degenerate; there are less than $2k^2$ degenerate quadruples in $A^4$, and at most $k^3$ non-degenerate ones. Now for a given quadruple, one can calculate the probability that it is also contained in $B^4$. This probability is $p^4$ if the quadruple is non-degenerate, and $p^2$ otherwise. Thus
\[\mathbb{E} M(B) \; = \;  \mathbb{E} \left(\sum s(x)^2\right) \; \leq \; 2k^2p^2 + k^3p^4.  \]
It follows immediately that
\begin{equation}\label{eq737} \mathbb{P} \left( M(B) \; \leq \; \frac{2k^2p^2 + k^3p^4}{1 - \epsilon/3}\right) \; \geq \; \epsilon/3.\end{equation}
Comparing \eqref{eq717} and \eqref{eq737}, we see that there is a specific choice of $B$ for which
\begin{equation}\label{eq757} \left| |B| - pk \right| \; \leq \; \epsilon p k/3 \qquad \mbox{and} \qquad M(B) \; \leq \; \frac{2k^2p^2 + k^3p^4}{1 - \epsilon/3}.\end{equation}
Now the Cauchy-Schwarz inequality gives
\[ M(B) \; = \; \sum_{x \in B\hat{+} B} s(x)^2 \; \geq \; \frac{\left(\sum s(x)\right)^2}{|B \hat{+} B|} \; = \; \frac{|B|^2(|B|-1)^2}{|B \hat{+} B|}.\] 
This, together with \eqref{eq757} and a small amount of computation, shows that $|B \hat{+} B| \geq (1 - \epsilon)k$. Moreover we have, from \eqref{eq757}, that \[ |B| \; \leq \; pk(1 + \epsilon/3) \; \leq \; \frac{3\log(1/\epsilon)}{\epsilon}\sqrt{k}.\]
This completes the proof.\endproof
\begin{corollaryenv}
\label{trickcor}
Let $k$ be a sufficiently large positive integer and let $A$ be a subset of an abelian group with $|A| = k$. Then $A$ contains a subset $C$, $|C| \geq k/8$, with $|C \hat{+} C| \geq 7|C|$.
\end{corollaryenv}
\proof Apply Proposition \ref{trickprop} with $\epsilon = 1/9$. This supplies us with a set $B$, $|B| \leq 100\sqrt{k}$, such that $|B \hat{+} B| \geq 8k/9$. Now simply take $C$ to be the union of $B$ and some arbitrary elements of $A$ so as to ensure that $|C|$ lies between $k/8$ and $8k/63$.\endproof \\[11pt]
\noindent\textbf{6. The rank of sets with small doubling.} In this section we obtain results concerning the Freiman dimension of sets $A \subseteq \Gamma$ with $|A| = k$ and $|A \hat{+} A| = m$. It is in this section that the difference between different groups $\Gamma$ comes to the fore. When $\Gamma = \mathbb{Z}_N$ we can relate $A$ to a subset of the integers, but this is not possible for $\Gamma = \mathbb{Z}_2^n$ and, by necessity, we must find a different approach.\\[11pt]
Let $A \subseteq \mathbb{Z}_N$ have $|A| = k$ and $|A \hat{+} A| = m$ and consider the image $\overline{A} \subseteq \mathbb{Z}$ of $A$ under the least positive residue or ``unfolding'' map $\mathbb{Z}_N \hookrightarrow \{1,\dots,N\}$. This clearly has $|\overline{A}| = k$ and $|\overline{A} \hat{+} \overline{A}| \leq 2m$. Thus we may apply 
the following result of Freiman \cite[page 24]{Fr}.
\begin{lemmaenv}[Freiman] \label{FL}Let $S \subseteq \mathbb{Z}$ have cardinality $k$ and Freiman dimension $r = r_{\mathbb{Q}}(S)$.  Then we have the inequality
\[ |S \hat{+} S| \; \geq \; r\left(k - \frac{r+1}{2}\right).\]
\end{lemmaenv}
\begin{lemmaenv}\label{lemma22} Suppose that $A \subseteq \mathbb{Z}_N$ has $|A| = k$ and $|A \hat{+} A| = m$, and let $\overline{A}$ be the image of $A$ under the unfolding map $\mathbb{Z}_N \hookrightarrow \mathbb{Z}$. Then we have the bound $r_{\mathbb{Q}}(\overline{A}) \leq 4m/k$.
\end{lemmaenv}
\proof Observe that $r = r(\overline{A}) \leq k - 1$, this being a general property shared by all sets of $k$ integers. Applying Lemma \ref{FL}, we have
\begin{equation}\boxeq 2m \; \geq \; |\overline{A} \hat{+} \overline{A}| \; \geq \; r\left(k - \frac{r+1}{2}\right) \; \geq \; \frac{1}{2}kr.\end{equation} 
\begin{propositionenv}\label{ZNbound}
We have the bounds
\begin{equation}\label{bb1} |S_k^m(\mathbb{Z}_N)| \; \leq \; N^{1 + \frac{4m}{k}}\left(\frac{2em}{k}\right)^k\exp\left(k^{31/32}\right)\end{equation} if $m \leq k^{31/30}/2$ and
\begin{equation}\label{bb2} |S_k^m(\mathbb{Z}_N)| \; \leq \; N^{1 + \frac{4m}{k}}k^{4k} \end{equation} whatever the value of $m$.
\end{propositionenv}
\proof It suffices to count the number of lifts $\overline{A}$ of sets $A \in S_k^m(\mathbb{Z}_N)$. Let $X_1,\dots,X_M$ be a complete set of representatives of the isomorphism classes of subsets $X$ of the integers with $|X| = k$ and $|X \hat{+} X| \leq 2m$. We always have $M \leq k^{4k}$ by Lemma \ref{lem12} but, when $m \leq k^{31/30}/2$, Proposition \ref{prop18a} allows us to do better and we have the bound $M \leq (2em/k)^k\exp(k^{31/32})$. Now each $\overline{A}$ is the image of some $X_i$ under some isomorphism $\psi : X_i \rightarrow \{1,\dots,N\}$, and by Lemma \ref{lem4} (i) there are at most $N^{r_{\mathbb{Q}}(X_i) + 1}$ such isomorphisms. The proposition follows immediately from these observations and Lemma \ref{lemma22}.\endproof\\[11pt]
We now turn our attention to subsets of $\mathbb{Z}_2^n$.
An important tool will be the following consequence of an inequality of Pl\"unnecke which was reproved by Ruzsa \cite{Ruz}. If $G$ is any abelian group and if $X \subseteq G$ then we write $kX$ for the $k$-fold sum $X + X + \dots + X$ and $kX - lX$ for the set $X + \dots + X - X - \dots - X$, where there are $k-1$ plus signs and $l$ difference signs.
\begin{theoremenv}[Pl\"unnecke -- Ruzsa]\label{thm23} Let $A \subseteq G$ be such that $|A + A| \leq C|A|$. Then for any positive integers $k$ and $l$ one has
\[ |kA - lA| \; \leq \; C^{k+l}|A|.\]
\end{theoremenv}
This allows us to get a bound on the Freiman dimension over $\mathbb{Z}_2$ of subsets of $\mathbb{Z}_2^n$ with small doubling.
\begin{lemmaenv}\label{mod2dim}
Suppose that $A \subseteq \mathbb{Z}_2^n$ has cardinality $k$ and $|A \hat{+} A| = m$. Then $r_{\mathbb{Z}_2}(A) + 1 \leq 4m\log k/k$.
\end{lemmaenv}
\proof
By Corollary \ref{cor1} we may assume (by subjecting $A$ to a Freiman isomorphism if necessary) that $A$ is not contained in any affine subspace of $\mathbb{Z}_2^n$ of dimension less than $r_{\mathbb{Z}_2}(A)$.\\[11pt]
Let $\prec$ be an arbitrary ordering on $\mathbb{Z}_2^n$ and let $h = \lfloor \log_e k\rfloor$. Suppose that $|A| = k$ and that $|A \hat{+} A| = m$, so that $|A + A| \leq m + 1$ (we are in characteristic 2!), and let $X \subseteq A$ be a subset of $A$ which is such that the sums $x_1 + x_2 + \dots + x_h$ ($x_1 \prec \dots \prec x_h$) are all distinct, and which is maximal with respect to this property. If $|X| \geq h$ then we have, using Theorem \ref{thm23},
\begin{eqnarray*}
\left(\frac{|X|}{h}\right)^h & \leq & \binom{|X|}{h}\\ & \leq & |hX| \\ & \leq & |hA| \\ & \leq & \left(\frac{m + 1}{k}\right)^hk.\end{eqnarray*}
It follows immediately from this that
\begin{equation}\label{eq11} |X| \; \leq \; 2mhk^{-1+1/h} \; \leq \; 4m\log k/k.\end{equation} 
Now the maximality of $X$ implies that if we set $Y = X \cup \{a\}$ for any $a \in A \setminus X$ then some two of the sums $y_1 + \dots + y_h$, $y'_1 + \dots + y'_h$ must be equal, where $y_1 \prec \dots \prec y_h$ and $y'_1 \prec \dots \prec y'_h$. It is clear that precisely one of the $y_i$'s or $y'_i$'s must equal $a$, and so $a \in hX - (h-1)X$. In particular, $A$ is contained in the affine subspace of $\mathbb{Z}_2^n$ spanned by $X$. The required inequality follows.\endproof\\[11pt]
It is interesting to compare this with the estimate
\[ r_{\mathbb{Z}_2}(A) + 1 \; \leq \; \log k + (m/k)^4 + 2\log(m/k).\] This bound, which follows from the work of Ruzsa \cite{RuzFrei} on Freiman's theorem in torsion groups, is stronger than Lemma \ref{mod2dim} when $m/k \ll (\log k)^{1/3}$.  Ruzsa and I, in work in progress, can improve this bound to something more like $\log k + 2(m/k)^2$. It would be interesting to know the true behaviour: perhaps $r \leq \log k + Cm/k$.
\begin{propositionenv}
\label{Z2bound}
We have the bounds
\begin{equation}\label{bb1a} |S_k^m(\mathbb{Z}_2^n)| \; \leq \; N^{\frac{4m\log k}{k}}\left(\frac{em}{k}\right)^k\exp\left(k^{31/32}\right)\end{equation} if $m \leq k^{31/30}$ and
\begin{equation}\label{bb2a} |S_k^m(\mathbb{Z}_2^n)| \; \leq \; N^{\frac{4m\log k}{k}}k^{4k} \end{equation} whetever the value of $m$.
\end{propositionenv}
\proof This is slightly more straightforward than Proposition \ref{ZNbound}. Let $X_1,\dots,X_M$ be a complete set of representatives of the isomorphism classes of sets $X \subseteq \mathbb{Z}_2^n$ with $|X| = k$ and $|X \hat{+} X| \leq m$. Proposition \ref{prop18a} gives bounds for $M$, and we know that any set $A \in S_k^m(\mathbb{Z}_2^n)$ is the image of some $X_i$ under an isomorphism $\psi : X_i \rightarrow \mathbb{Z}_2^n$. Lemma \ref{lem4}, (i) and Lemma \ref{mod2dim} allow us to bound the number of such isomorphisms.\endproof\\[11pt]
\noindent\textbf{7. Putting everything together.} In this section we use Propositions \ref{ZNbound} and \ref{Z2bound} to give upper bounds for the sums
\[ \sum_{m \geq 7k} |S_k^m(\Gamma)|2^{-m}\] when $\Gamma = \mathbb{Z}_N$ and $\mathbb{Z}_2^n$. An application of Proposition \ref{trickprop} will then allow us to prove Theorems \ref{mainthm1} and \ref{mainthm2}. We begin with Theorem \ref{mainthm1}. Let $k = \lfloor 20\log N\rfloor$. Using (\ref{bb1}) and (\ref{bb2}) we have
\begin{eqnarray*}
\sum_{m \geq 7k} |S_k^m(\mathbb{Z}_N)|2^{-m} & = & \sum_{m =7k}^{\lfloor k^{31/30}/2\rfloor} |S_k^m(\mathbb{Z}_N)|2^{-m} + \sum_{m =\lfloor k^{31/30}/2\rfloor}^{k(k-1)/2} |S_k^m(\mathbb{Z}_N)|2^{-m} \\ & \leq & \sum_{m = 7k}^{\lfloor k^{31/30}/2\rfloor} 2^{m\left(\left(\frac{4}{k} + \frac{1}{m}\right)\log N + \frac{k}{m}\log\left(\frac{2em}{k}\right) - 1 + o(1)\right)} + \sum_{m \geq \lfloor k^{31/30}/2\rfloor} 2^{m\left(\frac{4\log N}{k} - 1 + o(1)\right)}\end{eqnarray*}
Now $\log(2eC)/C - 1 \leq -0.2499$ when $C \geq 7$, and so both of the sums appearing here are bounded above by $N^{-2}$ for $N$ sufficiently large.\\[11pt]
Now suppose that some randomly selected set $A \subseteq \mathbb{Z}_N$ contains $X \hat{+} X$ with $|X| = k$. Then, by Proposition \ref{trickprop}, $A$ contains some set $Y \hat{+} Y$ where $|Y| \geq k/8$ and $|Y \hat{+} Y| \geq 7|Y|$. However if $k = \lceil 160 \log N\rceil$ then the expected number of such sets $Y$ that $A$ contains is at most $N^{-2}$, and in particular the probability that it contains even one is $o(1)$. It follows immediately that almost surely $A$ does not contain any $X \hat{+} X$ with $|X| = \lceil 160 \log N\rceil$, and this is the claim of Theorem \ref{mainthm1}.\endproof\\[11pt]
Similar computations apply for the group $\mathbb{Z}_2^n$, but now of course we use formulae (\ref{bb1a}) and (\ref{bb2a}) of Proposition \ref{Z2bound} instead. Let $k = \lfloor 11\log N\log\log N\rfloor$. One has
\begin{eqnarray*}
\sum_{m \geq 7k} |S_k^m(\mathbb{Z}_2^n)|2^{-m} & = & \sum_{m =7k}^{\lfloor k^{31/30}/2\rfloor} |S_k^m(\mathbb{Z}_2^n)|2^{-m} + \sum_{m =\lfloor k^{31/30}/2\rfloor}^{k(k-1)/2} |S_k^m(\mathbb{Z}_2^n)|2^{-m} \\ & \leq & \sum_{m = 7k}^{\lfloor k^{31/30}/2\rfloor} 2^{m\left(\frac{4\log k\log N}{k} + \frac{k}{m}\log\left(\frac{em}{k}\right) - 1 + o(1)\right)} + \sum_{m \geq \lfloor k^{31/30}/2\rfloor} 2^{m\left(\frac{4\log k\log N}{k} - 1 + o(1)\right)}\end{eqnarray*}
Now $\log(eC)/C - 1 \leq -0.39$ when $C \geq 7$, and so both sums appearing here are bounded above by $N^{-2}$ for $N$ sufficiently large. The rest of the argument proceeds exactly as before, and one gets Theorem \ref{mainthm2} (with $C_2 = 308$).\endproof\\[11pt]
\noindent\textbf{8. A lower bound for $\mathbb{Z}_2^n$.} In this section we prove the lower bound in Theorem \ref{mainthm2}. Throughout this section we write $\Gamma = \mathbb{Z}_2^n$ and $N = |\Gamma|$. 
\begin{propositionenv} \label{prop47}
\[\lim_{n \rightarrow \infty}\mathbb{P}_{\Gamma}(\omega \geq \textstyle\frac{1}{4}\displaystyle\log N\log\log N) \; = \; 1.\]
\end{propositionenv}
\proof Let $k = \lfloor \log n + \log\log n - 1\rfloor$ and consider all $k$-dimensional subspaces $H_1,\dots,H_M$ of $\Gamma$, where
\[ M \; = \; \frac{(2^n - 1)(2^n - 2) \dots (2^n - 2^{k-1})}{(2^k - 1)(2^k - 2) \dots (2^k - 2^{k-1})} \; \geq \; 2^{nk - k^2}.\] For a random set $A \subseteq \Gamma$ let $X_i$ denote the indicator function of the event $H_i \setminus \{0\}\subseteq A$, and let $X = X_1 + \dots + X_M$. We have
\begin{equation}\label{expX} \mathbb{E}X \; \geq \; 2^{nk -k^2 - 2^k}.\end{equation}
We will also estimate $\Var(X)$, and to do this it is necessary to prove a lemma.
\begin{lemmaenv}
The number of pairs $(i,j)$ for which $\dim (H_i \cap H_j) = l$, $0 \leq l \leq k$,
is at most $2^{n(2k - l)}$.
\end{lemmaenv}
\noindent\textbf{Proof of Lemma.} We will show that the quantity in question is exactly
\[ \frac{(2^n - 1)(2^n - 2)\dots(2^n - 2^{2k - l - 1})}{(2^k - 2^l)^2\dots(2^k - 2^{k-1})^2(2^l - 1)(2^l - 2) \dots (2^l - 2^{l-1})},\] which is easily seen to be at most $2^{n(2k - l)}$.
To make this count, first pick the intersection $L = H_i \cap H_j$ in any of
\[ \frac{(2^n - 1)(2^n - 2) \dots (2^n - 2^{l-1})}{(2^l - 1)(2^l - 2) \dots (2^l - 2^{l-1})} \] ways. Then pick vectors $u_1,\dots,u_{k-l},v_1,\dots,v_{k-l}$ so that $H_i = \Span(L,u_1,\dots,u_{k-l})$ and $H_j = \Span(L,v_1,\dots,v_{k-l})$. The subspace $L$ together with all the points $u_i,v_i$ spans a subspace of dimension $2k -l$, so the number of choices for the $u_i,v_i$ is
\[ (2^n - 2^l)(2^n - 2^{l+1})\dots(2^n - 2^{2k - l - 1}).\]
However, many of these choices will give the same subspaces $H_i,H_j$. If $H_i,H_j$ are fixed and have $H_i \cap H_j = L$ then there are $(2^k - 2^l)\dots(2^k - 2^{k-1})$ ways of choosing $u_1,\dots,u_{k-l}$ so that $H_i = \Span(L,u_1,\dots,u_{k-l})$, and the same number of ways of choosing $v_1,\dots,v_{k-l}$. This completes the proof of the lemma.\endproof\\[11pt]
Now the covariance $\Cov(X_i,X_j)$ is positive if and only if $H_i\setminus\{0\}$ intersects $H_j\setminus \{0\}$, which requires $\dim(H_i \cap H_j) \geq 1$. Thus we have, using the lemma,
\begin{eqnarray*}
\Var(X) & = & \sum_{i,j} \Cov(X_i,X_j) \\ & \leq & 2\sum_{l = 1}^k 2^{n(2k - l)}2^{-2^{k+1}}\left(2^{2^l} - 1\right) \\ & \leq & 2\sum_{l = 1}^k 2^{2kn - 2^{k+1} + 2^l - nl}.
\end{eqnarray*}
Now one can check that $2^l - nl \leq 2 - n$ for $1 \leq l \leq k$, and so
\[ \Var(X) \; \leq \; 8k2^{(2k - 1)n - 2^{k+1}}.\]
Now Chebyshev's inequality tells us that
\[ \mathbb{P}(X = 0) \; \leq \; \frac{\Var(X)}{(\mathbb{E}X)^2} \; \leq \; 8k2^{2k^2 - n},\]
which is $o(1)$. This means that with probability $1 - o(1)$ the set $A$ contains some $H_i \setminus \{0\}$. Note, however, that $H_i \setminus \{0\}$ is precisely $H_i \hat{+} H_i$.\endproof\\[11pt]
\noindent\textbf{9. Concluding remarks.} It would certainly be of interest to determine the clique number of a random Cayley sum graph on $\mathbb{Z}_N$ asymptotically (almost surely). There are several points at which are argument leading to Theorem \ref{mainthm1} was suboptimal. One was that we rather shied away from the issue of counting sets with extremely small doubling property, resorting instead to the trick of Proposition \ref{trickprop}. It is possible to say rather more about these using Freiman's deep theorem \cite{Fr}, which says that such sets are economically contained in a multidimensional arithmetic progression. It was also rather wasteful to treat subsets of $\mathbb{Z}_N$ by using the unfolding map $\mathbb{Z}_N \hookrightarrow \{1,\dots,N\}$ in such a na\"{\i}ve way. When $N$ is prime, it is possible to prove the following result:
\begin{propositionenv}\label{prop18} Suppose that $N$ is prime and that $A$ is a subset of $\mathbb{Z}_N$ with $|A \hat{+} A| \leq (\log N)^2/32$.
Then there are $\lambda \in \mathbb{Z}_N^{*},\mu \in \mathbb{Z}_N$ for which $\lambda A + \mu \subseteq \{0,\dots,\lfloor N/2\rfloor\}$.
\end{propositionenv}
By subjecting our set $A$ to such a homothety one can ensure that the unfolding map induces a Freiman isomorphism, so that the lifting $\overline{A}$ has $|\overline{A} \hat{+} \overline{A}| = m$.
Proposition \ref{prop18} may be called a \textit{rectification principle} in the spirit of \cite{BLR}.\\[11pt]
Using these results, and tidying up the approximations used in, for example, Lemma \ref{lemma22}, one can certainly reduce the number 160 drastically, and perhaps even to $(3 + \epsilon)$. To do this would be difficult, and there seem to be some even more substantial obstacles to pushing the constant down to $(2 + \epsilon)$, which seems likely to be the correct value. We may return to these issues in a future paper.\\[11pt]
Another issue deserving of remark arises from the result of Graham and Ringrose \cite{GR}, which implies that there is an absolute constant $c$ such that, for infinitely many primes $N$, the clique number $\omega(P_N)$ is at least $c\log N\log\log\log N$. However, using the Borel-Cantelli lemma together with the bounds obtained in \S 7 one can show that almost surely any sequence $(G_N)_{N}$ of random Cayley sum graphs on $\mathbb{Z}_N$, $|G_N| = N$, will have $\limsup_{N \rightarrow \infty} \omega(G_N)/\log N \leq 160$. Thus Paley graphs are (in a sense) not as good Ramsey graphs as random Cayley sum graphs.\\[11pt]
Let us close by mentioning a conjecture from \cite{Alon}.
\begin{conjectureenv}[Alon] There exists an absolute constant $b$ such that for every group $\Gamma$ on $n$ elements there is a Cayley graph on this group containing neither a complete subgraph nor an independent set on more than $b\log n$ vertices.\end{conjectureenv}
The results of this paper confirm this conjecture for cyclic groups but one should not expect the full conjecture, if it is true, to be easy to prove. Quite apart from the fact that our methods seem to break down hopelessly for all but a few very special groups, the result of Proposition \ref{prop47} shows that random methods alone will not suffice to attack Alon's conjecture.\\[11pt]
\noindent\textbf{10. Acknowledgements.} I would like to thank Imre Ruzsa and Seva Lev for their generous help in improving the exposition in this paper.

\noindent Ben Green\\
Trinity College, Cambridge, England.\\
email: bjg23@hermes.cam.ac.uk

     \end{document}